\documentclass [11pt]{amsart}
\usepackage {amsmath, amssymb, a4wide, epsfig}
\usepackage[latin1]{inputenc}
\date{December 22, 2004}

\author{Romain Dujardin}
\title{Structure properties of  laminar currents on $\pd$}

\newcommand{\cc}{\mathbb{C}}

\newcommand{\dd}{\mathbb{D}}

\newcommand{\e}{\varepsilon}
\newcommand{\cv}{\rightarrow}
\newcommand{\cvf}{\rightharpoonup}
\newcommand{\fr}{\partial}
\newcommand{\om}{\Omega}
\newcommand{\set}[1]{\left\{#1\right\}}
\newcommand{\norm}[1]{\left\Vert#1\right\Vert}
\newcommand{\abs}[1]{\left\vert#1\right\vert}
\newcommand{\finc}{\subset \subset}
\newcommand{\cd}{\cc^2}
\newcommand{\pd}{{\mathbb{P}^2}}
\newcommand{\pu}{{\mathbb{P}^1}}
\newcommand{\rest}[1]{ \arrowvert_{#1}}
\newcommand{\m}{{\bf M}}

\newcommand{\unsur}[1]{\frac{1}{#1}}
\newcommand{\el}{\mathcal{L}}

\newcommand{\qq}{\mathcal{Q}}

\newcommand{\ms}{{\rm m.s.}} 

\newcommand{\dft}{{\rm dft}}
\newcommand{\dint}[1]{\text{d-int}(#1)}

\DeclareMathOperator{\supp}{Supp}

\newtheorem{prop} {Proposition} [section]
\newtheorem{thm}[prop] {Theorem} 
\newtheorem{defi}[prop] {Definition}
\newtheorem{lem}[prop] {Lemma}
\newtheorem{cor}[prop]{Corollary}

\newtheorem{exam}[prop]{Example}

\newtheorem{rmk}[prop]{Remark}

\newenvironment{pf}{\noindent {\bf Proof:} }{\hfill $\square$\\}

\begin{document}

\begin{abstract}
We study the structure of a class of laminar closed positive currents on 
$\mathbb{CP}^2$, naturally appearing in birational dynamics. 
We prove such  a current
admits natural non intersecting {\em leaves}, that are closed under analytic continuation.
As a consequence it  can be seen as a  {\em foliation cycle}  a weak lamination.
\end{abstract}

\maketitle

\section{Introduction}

 Positive closed currents play an important role in higher dimensional holomorphic
 dynamics. Since the very beginnings of the theory, they have served as an elementary 
 bridge between the ambient complex geometry and the dynamics. Let us focus on the case of
 polynomial automorphisms of $\cd$ (see N. Sibony \cite{si} for a more thorough study and 
 bibliographical data). In case $f$ is hyperbolic on its nonwandering set, the laminar structure of 
 the Julia sets $J^+$ and $J^-$ is predicted by general Stable Manifold theory, whereas D. Ruelle and D. Sullivan proved in \cite{rs} the existence of foliation cycles (uniformly laminar currents) subordinate 
 to those laminations. E. Bedford and J. Smillie proved in \cite{bs1} that the Ruelle-Sullivan
currents coincide with the invariant ``Green" 
currents obtained by equidistributing preimages
of generic subvarieties. 

Laminar currents were introduced by E. Bedford, M. Lyubich and J. Smillie \cite{bls}
as analogues of the Ruelle-Sullivan foliation cycles in the general (non-uniformly hyperbolic) setting. They proved the invariant currents of polynomial automorphisms of $\cd$ are laminar and derived 
some dynamical consequences, among them the local product structure of the maximal 
entropy measure, as well as equidistribution of saddle periodic orbits. Also, the laminar structure of the Julia sets is required  in the notions of unstable critical points and external rays (see \cite{bs5, bs6}).\\
 
 Our purpose here is to continue the development of the general theory 
 of laminar currents with a view to new dynamical applications. 
 We proved in \cite{lamin} that laminar currents are abundant in rational dynamics on the 
 complex projective plane, by exhibiting a general criterion ensuring 
 laminarity of the limit of a sequence of $\mathbb{Q}$-divisors on $\pd$. We proved in \cite{isect} that 
 these currents are well behaved with respect to taking wedge products. 
 
 In this paper we give a structure theorem for this class of laminar currents, 
 which is also new
in the case of  plane polynomial automorphisms. In the article \cite{birat}, 
we apply these results to birational surface dynamics. \\
 
 Let us be more specific. Recall  a {\sl laminar current} in an open $\om\subset\cd$ is 
 a current ``filled" in the sense of measure, by compatible holomorphic disks (see below section \ref{sec_prel} for more details). The point is that the disks are not assumed to be
 properly embedded in $\om$. On the other hand a current is said to be {\sl uniformly laminar}
 if it is locally made up of currents of integration over disjoint complex {\sl submanifolds}. 
 
 There are examples showing that closed laminar currents may have somehow
 strange structure; this prevents us from studying {\sl general}
 laminar currents.
Nevertheless, the currents arising in some dynamical situations
(e.g. those constructed in \cite{bs5,lamin}) have an additional property: they are limits of 
sequences of
$\mathbb{Q}$-divisors $\unsur{d_{n}}[C_{n}]$ on $\pd$ with controlled geometry. We call such currents 
{\sl strongly approximable}. An important class of examples is provided by invariant currents of polynomial automorphisms, and more generally, birational maps. 
Our main result is the following: 

\begin{thm}\label{thm_main}
Let $T$ be a diffuse strongly approximable laminar current on $\pd$. Then
\begin{enumerate}
 \item If $\el\subset\om\subset \pd$ is an embedded lamination by Riemann surfaces, then 
 $T\rest{\el}$ is uniformly laminar (analytic continuation statement).
 \item Two disks subordinate to $T$ are compatible, i.e. their intersection is either empty or open in the disk topology (non self-intersection).
 \end{enumerate}
 \end{thm}
 
 Here, the notion of disk subordinate to $T$ is stronger than just appearing in the  decomposition 
 of $T$ as integral over a measured family of disks. It appears that the good notion to be considered is that of {\sl uniformly laminar current} subordinate to $T$ (see Definition
 \ref{def_sub} below).
 Note that in the case of non diffuse currents, that is currents giving mass to algebraic curves, if $\el$ is a curve, 
 the first item is a consequence of the Skoda-El Mir Extension
 Theorem  (see Demailly \cite{de}).
  A main issue in this theorem is that there is no regularity hypothesis  on the potentials, which seems important  in view of wide dynamical applications. On the other hand assuming the wedge product  $T\wedge T$ (is well defined and) vanishes, then the second item in the theorem is automatic --this was shown to hold in case  $T$ admits local continuous potentials in \cite{isect}.\\
 
We use this result to show, using  a construction of Meiyu Su \cite{su}, that disks subordinate to such 
$T$ form a lamination, in a weak sense, and $T$ induces an invariant transverse measure on this lamination. We believe this construction provides a useful language for treating problems related to 
strongly approximable currents. In particular this clarifies the question of  differing
{\sl representations} of a laminar current by measurable families of disks. 

As an application we prove in section \ref{sec_su} that if such a 
$T$ is extremal --which is common in dynamical situations-- then the transverse measure is ergodic.  
We also apply item 2. of the Theorem to prove (Theorem \ref{thm_harmo}) that the potential of a strongly approximable current is 
either harmonic or identically $-\infty$ on almost every  leaf. \\
 
 The precise outline of the paper is as follows. In {\S} \ref{sec_prel}
 we recall 
some basic notions  related to  laminar currents. In {\S} \ref{sec_cont}
 and 
\ref{sec_isect} we prove our main theorem, whereas the interpretation in terms
 of weak lamination structure is given in {\S}\ref{sec_su}. We also
 relate  invariant transverse measures on 
 the weak laminations and closed laminar currents dominated by $T$. In section 
 \ref{sec_pluripotential} we give some applications of our study, with
 some pluripotential theoretic flavour: we study the potential of  
strongly approximable laminar currents 
 along the leaves (Theorem \ref{thm_harmo}) and prove such currents
  decompose as  sums of two closed laminar currents, one not charging pluripolar sets, 
 and the other  with full mass on a pluripolar set (Theorem \ref{thm_decomp}).\\
 
\noindent{\bf Remark.} Much of this work was initiated as part of the author's ph.D thesis \cite{these}, 
 under the direction of N. Sibony. We would like to thank him once again for his  help and advice.
 
 %%%%%%%%%%%%%%%%%%%%%%%%%%%%%%%%%%%
 
 \section{Preliminaries on laminar currents}
 \label{sec_prel}
 
 We begin by recalling some definitions and preparatory results on
laminar currents, that will be useful to us in the sequel. Additional references are [Du1-3,BLS]. A good  reference on
positive closed currents is Demailly's survey article \cite{de}.

The first definitions are local so we consider an open subset 
$\om \subset \cd$, and $T$ a positive $(1,1)$ current in $\om$.
We let $\supp(T)$ denote the (closed) support of $T$, $\norm{T}$ the trace 
measure and $\m (T)$ the mass norm; $[V]$ denotes the integration current over
the subvariety $V$, possibly with boundary.  Also
$\dd$ denotes the unit disk in $\cc$. 

\begin{defi}\label{def_ul}  
$T$ is uniformly laminar if for every $x \in \supp(T)$ there exists 
 open sets $V\supset U\ni x$, with $V$
biholomorphic to the unit bidisk 
$\mathbb{D}^2$ so that in this coordinate chart
$T \rest{U}$ is the direct integral of integration currents over 
a measured family of disjoint graphs in $\mathbb{D}^2$, i.e. :

there exists
a measure $\lambda$ on $\set{0}\times\dd$,  and a family $(f_a)$ 
of holomorphic functions  $f_a: \mathbb{D} \cv \mathbb{D}$ such
that $f_a(0)=a$, 
the graphs $\Gamma_{f_a}$ of two different $f_a$'s are disjoint,  and
\begin{equation}\label{eq_ul}
T \rest{U} = \int_{\set{0}\times \dd}
[\Gamma_{f_a}\cap U ]~ d\lambda (a).
\end{equation}
\end{defi}

A family of disjoint horizontal graphs
in  $\mathbb{D}^2$ form a lamination. 
More precisely
the holonomy map is H{\"o}lder continuous
in this case: it is a corollary of the celebrated 
$\Lambda$-lemma on holomorphic motions \cite{mss}, since laminations
by graphs in $\dd^2$ and holomorphic motions in the unit disk are two
sides of the same object. Another useful consequence of the theory of
holomorphic motions is that a lamination by graphs of some vertically
compact $X \subset \dd^2$ has an extension to a neighborhood of $X$
--this is a special case of a deep theorem of Slodkowski's, see \cite{sl,do}. 
A uniformly laminar current
induces an invariant  transverse measure on its underlying
lamination. 

\begin{defi}\label{def_l}
$T$ is laminar in $\om$ if there exists a sequence of open subsets
$\om^i\subset\om$, such that   
$\norm{T}(\fr\om^i)=0$, together with an increasing sequence of
currents $(T^i)_{i\geq 0}$, $T^i$ uniformly laminar in $\om^i$,
converging to $T$. 

Equivalently (see \cite{bls}), $T$ is laminar in $\om$ if there exists a measured family ($\mathcal{A},\mu$) 
of  holomorphic disks $D_a\subset\om$, such that for every pair $(a,b)$, $D_a\cap D_b$ is either 
empty or open in the disk topology (compatibility condition), and 
\begin{equation}\label{eq_rep}
T=\int_\mathcal{A} [D_a]d\mu(a).
\end{equation}
\end{defi}

Notice that both representations of a laminar current as increasing limits or integral of disks 
are far from being unique, since they may be modified on sets of zero $\norm{T}$ 
measure. 
One aspect of the results in this paper (in particular the construction in section
 \ref{sec_su}) is to provide 
a {\sl natural} representation of $T$ as a foliated cycle on a
measured lamination, which is the ``maximum" of all possible
representations. Notice also that our laminar currents are called
``weakly laminar" in \cite{bls}.  

There is a useful alternate representation of $T$ as an integral of disks. 
Following \cite[equation (6.5)]{bls}, one may
reparametrize the  laminar representation (\ref{eq_rep}) to obtain an alternate representation
 of $T$ as an integral over a family of {\sl disjoint} 
disks. Indeed there exists a measured family $(\tilde{\mathcal{A}},\tilde\mu)$
of disjoint disks $(D_a)_{a\in \tilde{\mathcal{A}}}$, 
and for almost every $a\in \tilde{\mathcal{A}}$ a function $p_a$, nonnegative and a.e. equal to  a lower semicontinuous  function such that 
\begin{equation}
\label{eq_rep_disj}
T=\int_{\tilde{\mathcal{A}}} p_a[D_a]d\tilde\mu(a).
\end{equation}

A basic feature here is that there is some choice to be made in choosing the collection 
of disks that appear in this representation.
We now define a non ambiguous notion of disk subordinate to $T$. 
It has the advantage of
being independent of the representation of $T$ as a laminar current.
 
\begin{defi}\label{def_sub}
A holomorphic disk $D$ is subordinate to $T$ if there exists in some $\om'\subset\om$
a uniformly laminar current 
$S\leq T$ with positive mass, such that $D\subset \supp(S)$ 
lies inside a leaf of the lamination induced by $S$.
\end{defi}

This definition is motivated by the fact that the usual ordering 
on positive currents is compatible with the laminar structure. This is the content of 
the next proposition, which is implicit in \cite[{\S}6]{bls} 
(see \cite{these} or \cite{ca} for a precise proof). 

Remark that, unlike the disks of representation 
(\ref{eq_rep}), disks subordinate to $T$ may 
intersect.  Such examples are provided by sums of 
uniformly laminar currents with transversals of zero area, 
see e.g. \cite[Example 2.2]{isect}.

\begin{prop}\label{prop_sub}
Let $T_1\leq T_2$ be laminar currents in $\om$. Assume $T_2$ has the 
representation $T_2=\int_{\mathcal{A}} p_a[D_a]d\mu(a)$. Then $T_1$ and $T_2$ have compatible 
representations in the sense that $T_1$ may be written as 
$T_1=\int_{\mathcal{A}} q_a[D_a]d\mu(a)$, with $q_a\leq p_a$ almost everywhere.
\end{prop}

The currents we will be interested in  in this article have
a crucial  additional  property: there is an explicit bound on the
``residual mass''   $\m(T-T^i)$.  We call such currents {\sl strongly
approximable}.  We need first introduce a few concepts.

Let us consider a sequence of (one dimensional)
analytic sets $[C_n]$ defined in some
neighborhood of $\overline\om$, with area $d_n$, and such that
$d_n^{-1}[C_n]\cvf T$. The disks of $T$ are to be obtained as cluster
values of sequences of graphs for some linear projection $\pi$. Let
$L$ be a complex line, transverse to the direction of the projection $\pi$, 
and  a subdivision $\qq$ of $L$ into squares of size $r$. For $Q\in
\qq$, we say that a
connected component $\Gamma$ of $\pi^{-1}(Q)\cap C_n$ is {\sl good} if $\pi:
\Gamma\cv Q$ is a homeomorphism, and the area of $\Gamma$ is bounded
by some universal  constant, {\sl bad} if not. The assumption on ${\rm
  area}(\Gamma)$  
ensures that families of good components form normal families.

\begin{defi}\label{def_sa} $T$ 
is a strongly approximable laminar current if there exist
a sequence $(C_n)$ of analytic subsets of some neighborhood $\mathcal{N}$ of
$\overline\om$, with
$d_n^{-1}[C_n]\cvf T$, 
 at least two distinct linear projections $\pi_j$, and a constant $C$
such that if $\qq$ is any subdivision of the projection basis $L$ into
squares of size $r$, and $T_{\qq,n}$ denotes the current 
made up of good components of $d_n^{-1}[C_n]$ in $\mathcal{N}$, one has the
following estimate in $\om$
\begin{equation}\label{eq_goodcp}
\left \langle d_n^{-1}[C_n]- T_{\qq,n},
\pi_j^*(\omega\rest{L})\right\rangle \leq Cr^2,
\end{equation}
where $\omega\rest{L}$ is the restriction of the ambient K{\"a}hler form
to the complex line $L$. \\

We say $T$ is strongly approximable in $\pd$ if the $C_n$ are plane algebraic curves 
and assumption  (\ref{eq_goodcp}) on good components holds for a
generic linear projection $\pd\backslash \set{p}\cv\pu$.
\end{defi}

Here ``generic" means ``for $p$ outside a countable union of Zariski closed subsets.
An important consequence of the definition is that such currents are closed.
The definition of strongly approximable currents (locally or on $\pd$), 
though seeming rather inelegant, is designed to fit with the constructions  
in  \cite{lamin} and \cite{bs5}, where it is of course satisfied. 

\begin{rmk}\normalfont \label{rmk_global}
We did deliberately state a local definition, for we believe this is the good setting 
for future applications. However, item 1. (analytic contination) 
of the main theorem \ref{thm_main}, even though it is
 a local result,  requires global hypotheses.
More precisely the important fact is the following: we need the approximating curves $C_n$
to have controlled number of intersection points with the fibers $\pi^{-1}(x)$, which is only 
known to hold in global situations. It would be interesting to extend it to the purely 
 local case, and 
more generally, using only the mass estimate (\ref{eq_mass})
below.  De Thelin \cite{dt} has a local approach to 
approximation of laminar currents, nevertheless we do not know if the crucial estimate 
(\ref{eq_mass}) is true in his case. 

In the remainder of the paper, we will have the occasion to deal with results that require
the global hypothesis, and those that do not, we will then
respectively speak of currents strongly 
 approximable in $\pd$ or in some $\om$.
\end{rmk}

Equation (\ref{eq_goodcp}) can be turned into a real mass estimate when
combining both projections. For a proof of the next proposition, see
\cite[Proposition 4.4]{isect}. 

\begin{prop}\label{prop_mass}
Let $T$ be a strongly approximable laminar current in $\om$. Fix $\om'\finc\om$, 
and $\pi_1$,
$\pi_2$ projections satisfying 
definition \ref{def_sa}. Then for any  
subdivisions   $\mathcal{S}_1$, $\mathcal{S}_2$ of the respective projection
bases into squares of size $r$, if 
$$\qq = \set{\pi_1^{-1}(s_1)\cap\pi_2^{-1}(s_2), (s_1, s_2) \in
  \mathcal{S}_1\times \mathcal{S}_2}$$ denotes the associated
subdivision of $\om$ into affine cubes of size $r$, there exists a
current  $T_\qq \leq T$ in a neighborhood of $\overline{\om}$, 
uniformly laminar in each $Q\in \qq$, 
and satisfying the estimate 
\begin{equation}\label{eq_mass}
\m (T-T_\qq)\leq C r^2,
\end{equation}
in $\om'$, with $C$ independent of $r$.
\end{prop}

%%%%%%%%%%%%%%%%%%%%%%%%%%%%

\section{The defect function and analytic continuation}\label{sec_cont}

The aim of this section is to prove the first part of Theorem \ref{thm_main}. We introduce a
notion of {\sl defect} of  a laminar current with respect to a projection, analogous to
the Ahlfors-Nevanlinna defect for entire functions : the defect measures the amount of good 
components in the slice mass of the current in some fiber. We already stressed in remark \ref{rmk_global} that this section uses global arguments; for ease of reading we consider the case of plane algebraic curves, nevertheless the cases of horizontal-like curves in the bidisk (see \cite{these}), or curves on an algebraic surface are similar --what is needed is the approximating curves to have 
a controlled  number of intersection points with the fibers of the projection.
 
The way to the proof is  quite 
simple, but precise formulation requires some care, and we apologize by advance for 
possible stylistic heaviness.\\

We begin by recalling the statement of the analytic continuation theorem. If $T$ is a laminar current 
in $\om$ there is a representation (\ref{eq_rep}) of $T$ as an integral of compatible disks
$$T=\int_\mathcal{A} [D_a]d\mu(a).$$ If $\el\subset\om'\subset\om$ is an embedded lamination, 
we define the restriction $T\rest{\el}$ as 
\begin{equation}\label{eq_rest}
T\rest{\el}= \int_{\set{a\in\mathcal{A}, ~D_a\subset L\in \el}} [D_a]d\mu(a),
\end{equation}
 where  the notation $L\in\el$ means ``$L$ is a leaf of $\el$". The current
 $T\rest{\el}$ is laminar in $\om'$.

\begin{thm}\label{thm_anal_cont}
Let $T$ be a strongly approximable laminar current on $\pd$. Assume $\el$ is an embedded 
lamination in $\om\subset\pd$, then $T\rest{\el}$ is uniformly laminar in $\om$.
\end{thm}

It may not seem so clear why this is an analytic continuation result. 
Recall the alternate representation
(\ref{eq_rep_disj})
$$
T=\int_{\tilde{\mathcal{A}}} p_a[D_a]d\tilde\mu(a).$$
The functions $p_a$ need not be locally constant, even if $T$ is closed: see  Demailly's example \cite[example 2.3]{lamin}. The assertion of the theorem is that in case $T$ is strongly 
approximable, the functions $p_a$ are {\sl globally} constant 
along disks subordinate to $T$ 
(in the sense of definition \ref{def_sub}).\\

The scheme of the proof is quite natural. The approximating curves form branched 
coverings (global assumption) of growing degree and branching over a line in $\pd$. Slicing Theory provides a way to ``count" the number of points in fibers as {\sl measures} on the fibers 
and we try to construct as many local sections as possible, matching over overlapping disks in the basis --the sections of the covering being uniformly laminar currents subordinate to $T$.\\

For this section let us fix 
 a sequence of curves $C_n\subset\pd$ of degree $d_n$
satisfying definition \ref{def_sa}, $T=\lim d_n^{-1}[C_n]$, and fix a linear projection $\pi_{p}:\pd\backslash\set{p}\cv\pu$, 
such that the Lelong number $\nu(T,p)$ vanishes. We also assume $p\notin C_n$ for every $n$.
For almost every line through $p$, one may define the slice 
$T\rest{L}$ which is a probability measure on $L$, moreover $T\wedge [L]$ is well defined and 
$T\wedge [L]=T\rest{L}$.
In this case, Siu's stability property of Lelong numbers by slicing
(see Demailly \cite{de}) shows that for almost every $L$,  
$T\rest{L}$ gives no mass to $\set{p}$. We also choose $p$ such that the set of vertical disks 
for the projection $\pi_{p}$ in the laminar decomposition of $T$ has measure zero. \\

We are given a lamination in $\om\subset\pd$, and we want to prove the restriction 
$T\rest{\el}$ 
is uniformly laminar. Since the problem is local (on $\el$), 
we may assume $\el$ is made up of graphs 
over the unit square $Q_0\subset\cc$ for some projection $\pi_p$ satisfying the requirements above; moreover we may assume $\el$ is vertically compact. So for
 now we denote $\pi_p$ by $\pi$ and restrict the problem to 
$\pi^{-1}(Q_0)\simeq Q_0\times \cc$. We consider the following 3 sequences of overlapping subdivisions of $\cc$, where $\qq$ is the standard subdivision 
(tesselation) of $\cc$ into
 translates of $Q_0$ and $r_k\cv 0$ --say $r_k=2^{-k}$--
$$\qq_k^0= r_k\qq,~ \qq_k^1= r_k\qq+
(\frac{r_k}{3}+\frac{2ir_k}{3}),~   \qq_k^2= r_k\qq+
(\frac{2r_k}{3}+\frac{ir_k}{3});$$ 
these subdivisions induce subdivisions of $Q_0$ that form a neighborhood basis of $Q_0$. For each $Q\in \qq^j_k$, we let $\mathcal{G}(Q,n)$  be the family of good components of $C_n$ 
over $Q$, and $$T_{Q,n}=\unsur{d_n}\!
 \sum_{ \Gamma\in \mathcal{G}(Q,n)}\![\Gamma];$$ if $\qq$ is one of the subdivisions $\qq^j_k$ 
 we also use the notation $T_{\qq,n}=\sum_{Q\in\qq} T_{Q,n}$. Also, for $Q\in \qq^j_k$, there is 
 a subsequence, still denoted by $n$, such that $T_{Q,n}\cvf T_Q$, where $T_Q$ is a uniformly laminar current in $Q\times \cc$: see e.g. \cite[proposition 3.4]{lamin} 
 --recall that by definition, good components form normal 
 families. We perform a diagonal extraction so that for {every} $j,k$ and every 
  $Q\in\qq^j_k$, 
 one has the convergence $T_{Q,n}\cv T_Q$. Let  $T_{\qq^j_k}= \sum 
_{Q\in\qq^j_k} T_Q$; from (\ref{eq_goodcp}) one deduces the  estimate
\begin{equation}\label{eq_goodcp2}
\langle T-T_{\qq^j_k}, \pi^*(idz\wedge d\overline{z})\rangle
\leq Cr_k^2. 
\end{equation} 
\\
 
 The currents $T_{Q,n}$, $T_Q$ have an important property of {\sl invariance of the slice mass}, 
 that  we now describe. For a vertical fiber $F_x=\pi^{-1}(x)$ one has
  $$T_{Q,n}\rest{F_x}=T_{Q,n}\wedge [F_x]=\sum_{y\in F_x\cap \mathcal{G}(Q,n)} \!\delta_y$$
(there are no multiplicities  because the current is made of good components). The mass of the 
slice measures is constant and denoted by $\ms(T_{Q,n})$. The uniformly laminar currents $T_Q$ do have the same property, and $T_Q\wedge [F_x]$ is the image on $F_x$ of the underlying transverse
measure of $T_Q$.

Now we claim that $\ms(T_{Q,n})\cv\ms(T_Q)$.  Convergence on compact subsets implies 
$\liminf \ms(T_{Q,n})\geq\ms(T_Q)$ and we need to check the other inequality. We know that 
$d_n^{-1}[C_n]\cvf T$ in $Q\times \cc$, so for almost every fiber $F_x$,
$d_n^{-1}[C_n]\wedge[F_x]\cvf T\wedge[F_x]$, and, as noted before, the hypothesis $\nu(T,p)=0$ implies $ T\wedge[F_x]$ 
is a probability measure on $F_x$. Now write $d_n^{-1}[C_n]=T_{Q,n}+R_{Q,n}$, and 
assume $R_{Q,n_{\ell}}\cvf R_Q$ so that  $T=T_Q+R_Q$. As $T$ is a current on $\pd$, 
$\ms(T)=1$ is well defined, and so is the case for $R_Q$. We conclude using the 
inequality $\liminf \ms(R_{Q,n_{\ell}})\geq\ms(R_Q)$.

\begin{defi}\label{def_dft}
For $Q\in\qq^j_k$,  $k\geq 0$, $0\leq j \leq 2$, one defines the defect of $Q$ by
$$\dft(Q)=1-\ms(T_Q).$$
 \end{defi}
 
 The reference to $T$ is implicit here.
 Since $\ms(T_{Q,n})\cv\ms(T_Q)$, the defect is the asymptotic 
 proportion of bad components over $Q$. One has the following properties: 
 
 \begin{prop}\label{prop_dft}~

\begin{itemize}
\item[\sl i.] $k,j$ being fixed, $\displaystyle{\sum_{Q\in\qq^j_k}
\dft(Q)\leq C}$;
\item[\sl ii.] if $Q'\subset Q$ then $\dft(Q')\leq \dft(Q)$.
\end{itemize}
\end{prop}

\begin{pf} for $Q\in\qq^j_k$ one has the estimate
 $$\lim_{n\cv\infty} \langle d_n^{-1}[C_n]-T_{Q,n},\mathbf{1}_{Q\times\cc} \pi^*
(id z\wedge d\overline z)\rangle = \dft(Q) r_k^2,$$ and  {\sl i.} is a consequence of 
estimate (\ref{eq_goodcp2}).

To prove the second point, note that good components over $Q$ are good components 
over $Q'$, so $\ms(T_{Q,n})\geq \ms(T_{Q',n})$, and let $n\cv\infty$.
\end{pf}

\begin{defi}\label{def_dft_ponct}
For  $x\in Q_0$, we let $\dft(x)=\lim \dft(Q_p)$, where $(Q_p)$ is any decreasing sequence of squares
such that $\set{x}=\cap_p Q_p$.

The fiber $\set{x}\times \cc$ is regular if $\dft(x)=0$, singular (resp. $\e$-singular) otherwise (resp. if $\dft(x)\geq \e$). 
 \end{defi}
 
 One easily checks the definition of $\dft(x)$ is non ambiguous using property {\sl ii.} of the preceding proposition.
 
 \begin{prop}
 There are at most countably many singular fibers; moreover 
  $$\sum_{x\in Q_0}  \dft(x)\leq 3C.$$
 \end{prop}
 
 \begin{pf} $j$ being fixed, the number of squares $Q\in\qq^j_k$ where $\dft(Q)\geq \e$ is 
 less than $C/\e$, which implies the bound $C/\e$ on the number of $\e$-singular
 fibers. The result follows.
 \end{pf}
 
 \begin{rmk}\normalfont
 We have no result on the structure of singular fibers. If a diffuse strongly approximable
 current has Lelong Number $\geq \e$ at some point $p$, generic slices through $p$ 
 have a Dirac mass $\e$ at $p$ and the fiber is $\e$-singular. 
On the other hand the invariant currents associated to polynomial automorphisms 
of $\cd$ have no singular fibers. Here is a rough argument: 
take some complex line $L$
in $\cd$ and iterate $L$ backwards. It is known that 
the iterates $d^{-n}[f^{-n}(L)]$ 
converge to the stable current $T^+$ which is laminar; since the wedge 
product $T^+\wedge T^-$ is geometric (i.e. described by intersection of disks,
see \cite{bls} or \cite{isect}), this implies 
$f^{-n}(L)$ intersects many disks of $T^-$ and so does $L$. In particular $L$ is not a 
singular fiber associated to $T^-$.
\end{rmk}

The following proposition is the basic link between defect and analytic continuation.

\begin{prop}\label{prop_prol}
Let $Q,Q'$ be two squares such that $Q\cap Q'\neq\emptyset$, 
$\dft(Q)\leq \alpha$, $\dft(Q')\leq \alpha'$, with $\alpha+\alpha'<1$. 
Then there exists a uniformly laminar current $T_{Q\cup Q'}$
in $(Q\cup Q')\times \cc$, such that $$T_{Q\cup Q'}\rest{Q}\leq T_Q,~ 
T_{Q\cup Q'}\rest{Q'}\leq T_{Q'}$$ and  $$\ms (T_{Q\cup Q'})\geq 1-\alpha
-\alpha'.$$
\end{prop}
 
\begin{pf} as $n\cv\infty$, 
 $\lim \ms(T_{Q,n}) \geq  1-\alpha$, and  $\ms(T_{Q,n})=
\unsur{d_n} \# \mathcal{G}(Q,n)$ is the number of good components over $Q$. So for 
$n\geq n(\e)$  there are at least 
 $(1-\alpha-\e)d_n$ (resp. $(1-\alpha'-\e)d_n$) good components over $Q$ (resp. $Q'$). 
 As the total number of components over $Q\cap Q'$ is bounded by $d_n$, at least
 $(1-\alpha-\alpha'-2\e)d_n$ components match over $Q\cap Q'$ for $n$ large, 
 giving rise to as many global good components over $Q\cup Q'$. Extracting a convergent subsequence of the sequence  
  $$\unsur{d_n}\sum_{\Gamma\in\mathcal{G}(Q\cup Q',n)}\! 
[\Gamma]$$ gives the desired $T_{Q\cup Q'}$.
\end{pf}

We inductively use this proposition to construct analytic continuation of disks along 
paths in $Q_0$. We assume all paths are continuous.
 
 \begin{defi}\label{def_prolpath}
 Let $\gamma:[0,1]\cv Q_0$ be an injective path. $T$ is said to have almost analytic continuation property up to $\e$ along $\gamma$ if there exists an open set $V_{\e}\supset
\gamma$ and a uniformly laminar current $T_{V_{\e}}\leq T$ made up of graphs over $V_{\e}$, and such that $\ms(T_{V_{\e}})\geq 1-\e$.

$T$ has the analytic continuation property along $\gamma$ if it has the almost analytic continuation property up to $\e$ for all $\e>0$.
\end{defi}

Proposition \ref{prop_prol} has the following corollary. The proof is left to the reader.

\begin{cor}\label{cor_prolpath}
Let $\gamma:[0,1]\cv Q_0$ be an injective path. Suppose  there exists $\e>0$, and 
a covering of $\gamma$ by a family of squares $\mathcal{F}_{\e}$ satisfying 
$$\sum_{Q\in\mathcal{F}_\e} \dft(Q)\leq \e.$$
Then $T$ has almost analytic continuation property up to $\e$ along $\gamma$.
\end{cor}

By definition the {\sl total defect} of $\gamma$ is the lower bound of the sums 
$\sum_{Q\in\mathcal{F}} \dft(Q)$ for families $\mathcal{F}$ of squares in $\qq^j_k$ covering $\gamma$. 
The next proposition is the crucial technical point in the proof of the theorem. 

 \begin{prop}\label{prop_percol}
 Let $x_1$ and $x_2$ be two points in $Q_0$. For every $\e>0$ there exists
 a path $\gamma_{\e}$ such that $T$ has almost continuation 
 property along $\gamma_{\e}$ up to $(\dft(x_1)+\dft(x_2)+\e)$. 
 \end{prop}
 
 In particular if $x_1$ and $x_2$ are regular,  $T$ has the
  $\e$-almost continuation property. For ease of reading we use the following notations:
  $a\approx b$ means $c^{-1}a\leq b\leq ca$ and $a\lesssim b$ means $a\leq cb$, with $c$   a contant independent of $r$ (size of the squares). The idea of the proof is to construct 
  enough essentially disjoint paths joining $x_1$ and $x_2$ in the subdivisions and apply 
  proposition \ref{prop_dft}.
  
  \begin{pf} assume first $x_1$ and $x_2$ lie on the same horizontal line. 
 Consider ``big" squares $Q_1\ni x_1$ and $Q_2\ni x_2$ of size $\approx \sqrt{r}$, and 
 join $Q_1$ and $Q_2$ by $N\approx 1/\sqrt{r}$ disjoint horizontal paths $\gamma_i$ 
with mutual distance $\geq 10r$. Each path $\gamma_i$ is covered by a family of 
``small" squares of size $r$, $\mathcal{F}_i\subset \qq^0\cup\qq^1\cup \qq^2$.
The families $\mathcal{F}_i$ are disjoint.
 Complete the paths $\gamma_i$ by adding affine pieces so that the paths join $x_1$ and $x_2$.
 
 We now evaluate the total defect of the family of paths, using proposition \ref{prop_dft}
 $$\sum_{i=1}^N\dft(\gamma_i)\leq \sum_{i=1}^N \left(\dft(Q_1) + \dft(Q_2)+
  \!\sum_{Q\in\mathcal{F}_i}\!\dft(Q)\right) \leq N(\dft(Q_1) +
  \dft(Q_2)) + 3C.$$ 
 As $N\approx 1/\sqrt{r}$, the average defect is 
 $$\unsur{N}\sum_{i=1}^N\dft(\gamma_i)
\leq \dft(Q_1) + \dft(Q_2) +c \sqrt{r},$$
 so at least one of the paths has total defect $\leq \dft(Q_1) + \dft(Q_2) +c\sqrt{r}$.
 
 In the general case consider an affine isometry $h$ such  that $h(x_1)$ and $h(x_2)$ 
 lie on the same horizontal line, and remark that if  $Q$ is one of the squares 
  of the preceding construction $h(Q)$ is included in a square of size at most twice that of $Q$. 
 \end{pf}
 
We will prove  theorem \ref{thm_anal_cont} through the following reformulation, 
which is of independent interest. 

\begin{prop}\label{prop_reformul}
Let $U\subset Q_0$ be a connected open subset.
Let $S$ be a uniformly laminar current with vertically compact support in $U\times \cc$, made up of graphs over $U$, and such that $S\leq T$ in $U_1\times\cc$, for some open 
$U_1\subset U$.

Then $S\leq T$ in $U\times \cc$.
\end{prop}
 
The conclusion of the proposition is that the relation $S\leq T$ propagates to the domain of definition of $S$; this is a continuation result in terms of disks subordinate to $T$.\\

We first prove that this proposition implies the theorem. Recall we localized the problem so that $\el$ is a vertically compact lamination whose leaves are graphs $\Delta_\alpha$ 
over $Q_0$. Then 
one has 
$$T\rest{\el} = \int_X p_{\alpha} [\Delta_{\alpha}] d\mu_X(\alpha)$$ where $\mu_X$ is a
positive measure on the global transversal $X$ and for every $\alpha$, $p_\alpha$ is 
a.e. equal to a  lower semicontinuous 
nonnegative function. We have to prove the functions $p_\alpha$ are constant for a.e. $\alpha$. The idea is as follows: if $p_\alpha$ takes the value $p_0$ at some point, then the preceding proposition forces $p_\alpha\geq p_0$ on $\Delta_\alpha$. 

Indeed consider the measurable function 
$$X\ni\alpha\longmapsto \inf(p_\alpha),$$ 
where $\inf$ denotes essential infimum. Then for $\e>0$,
$$\Delta_\alpha^\e=\set{ x\in \Delta_\alpha, 
~p_\alpha(x)>\inf (p_\alpha)+\e}$$ is either empty or an open subset up to a set of zero area. We prove $\Delta_\alpha^\e$ has area zero for a.e. $\alpha$.

Assume the contrary. There exists $Y\subset X$ of positive transverse measure such that 
$\Delta_\alpha^\e$ has positive area. Hence the current  
$$T_\e=\int_Y\mathbf{1}_{\Delta_\alpha^\e}
p_\alpha[\Delta_\alpha] d\mu_X(\alpha)$$ is a laminar current of positive mass. By a monotone convergence argument and the Fubini theorem (see e.g. \cite[Proposition 6.2]{bls}) there exists a square $R$ and a set $Y_R$ of positive transverse
 measure such that for
$\alpha\in Y_R$, $\pi^{-1}(R)\cap \Delta_\alpha\subset \Delta_\alpha^\e$ up to a set of 
zero area. In particular 
$$0<S=\left(\int_{Y_R}(\inf (p_\alpha)+\e) d\mu_X(\alpha)\right) \rest{R\times\cc}
\leq T_\e \rest{R\times\cc}\leq T\rest{R\times\cc}.$$
Since $S$ is uniformly laminar, by proposition \ref{prop_reformul}, this relation propagates to $Q_0\times \cc$, contradicting the definition of $\inf p_\alpha$. \hfill $\square$\\

\noindent{\bf Proof of proposition \ref{prop_reformul}:} without loss of generality, assume
$U_1$ is a disk. We use the representation (\ref{eq_rep_disj}) over families of disjoint disks 
$T=\int_{\mathcal{A}} p_a[D_a]d\mu(a)$. Using proposition \ref{prop_sub} one gets 
$$S=\int_{\mathcal{A}\rest{U_1}} q_a [D_a] d\mu(a)$$ where $\mathcal{A}\rest{U_1}$ is the set of restrictions to $U_1\times \cc$ of the disks of $\mathcal{A}$, and $q_a$ is a constant on every $D_a\in\mathcal{A}\rest{U_1}$ since $S$ is uniformly laminar; moreover 
$q_a\leq p_a$ a.e. For a square $Q$ let $\mathcal{A}_Q\subset \mathcal{A}\rest{Q}$ be the set of disks that are graphs over $Q$. The uniformly laminar
currents $T_Q$ previously considered have the form
$$ T_Q=\int_{\mathcal{A}_Q} p_{a,Q}[D_a]d\mu(a)$$
where $p_{a,Q}$ is a constant function.

\begin{lem}\label{lem_pp}
For almost every $x\in U_1\times \cc$ there exists a decreasing sequence of squares
$Q_p$, $\cap_{p\geq 0} Q_p=\set{x}$, and for every $p$ a uniformly laminar current
$S_{Q_p}\leq S$, in $Q_p\times \cc$,
such that  $S_{Q_p}\leq T_{Q_p}$ et $\ms (S_{Q_p})\cv m=\ms(S)$ as $p\cv\infty$.
\end{lem}

\begin{pf} 
Recall that if $\qq_k$ is one of the 3 sequences of subdivisions $\qq^j_k$,
the sequence $T_{\qq_k}=\sum_{Q\in\qq_k}T_Q$ increases to $T$. Let 
$$S_{\qq_k}=\sum_{Q\in\qq_k} \int_{\mathcal{A}_Q}
\inf(q_a, p_{a,Q}) [D_a]d\mu(a)= \int_\mathcal{A} \inf \big(q_a, 
\sum_{Q\in\qq_k}\mathbf{1}_{Q\times\cc}~ p_{a,Q}\big) 
[D_a]d\mu(a).$$
The current $S_{\qq_k}$ is uniformly laminar in each $Q\times \cc$, $Q\in\qq_q$, and since 
$$p_{a,\qq}= \sum_{Q\in\qq_k}\mathbf{1}_{Q\times\cc}~ p_{a,Q}$$ increases $\norm{T}$ a.e. to $p_a\geq q_a$, one gets $\inf(q_a, p_{a,
\qq_k})\nearrow q_a$ and
the sequence of currents $S_{\qq_k}$ increases to $S$ as $k\cv\infty$.
From this one easily deduces the conclusion of the lemma
\end{pf}
 
 We continue with the proof of proposition \ref{prop_reformul}. The basic idea is to 
 transport the relation $S_Q\leq T_Q$ by using analytic continuation along paths.
 Fix $\e>0$ and $x_0\in U_1$ such that the conclusion of the lemma is satisfied and 
 $\dft(x_0)=0$. For $x_0\in Q\in \qq_k$, $k$ large enough, one has a uniformly laminar 
 $S_Q\leq T_Q$, such that $\ms(S_Q)\geq \ms(S)-\e=m-\e$.
 
 On the other hand by proposition \ref{prop_percol} for every $x_1\in U$ there exists a 
 path $\gamma_\e$ joining $x_0$ and $x_1$, such that $T$ has $(\dft(x_1)+\e)$-almost analytic continuation along $\gamma_\e$, i.e. there exists $V\supset\gamma_\e$, and $T_V\leq T$ uniformly laminar, such that $\ms(T_V)\geq 1-\dft(x_1)-\e$. Applying proposition \ref{prop_prol}  to $Q,V$, and the square $Q_1\in \qq_k$ containing $x_1$ yields the existence of a uniformly laminar current $T_{Q,Q_1}$, simultaneously 
subordinate to $T_Q$, $T_V$ and $T_{Q_1}$, 
with  $$\ms(T_{Q,Q_1})\geq 1-\dft(Q)-\dft(Q_1)-\e -\dft(x_1)\geq 
1-\dft(Q)-2\dft(Q_1)-\e.$$ By construction the graphs of   $T_{Q,Q_1}$ over $Q_1$ are the analytic continuations along $\gamma_\e$  of those over $Q$.\\

We then prove the sum of $T_{Q,Q_1}$, with varying $Q_1$, approximate $T$ in $U\times \cc$: 
\begin{align*}
\left \langle T- \sum_{Q_1\in\qq_k} T_{Q,Q_1},
  \mathbf{1}_{U\times\cc} ~\pi^* idz\wedge d\bar z
\right\rangle &\leq
  \sum_{Q_1\in\qq_k} (\dft(Q)+\e+2\dft(Q_1)){\rm area}(Q_1)\\ 
&\leq (\dft(Q) +\e) + 2 \left(\sum_{Q_1\in\qq_k} \dft(Q_1)\right) r_k^2
\end{align*}
and the right hand side is less than $3\e$ if $k$ is large and $Q$ small.

We now claim there exists for all $Q_1$ a current $S_{Q,Q_1}$ such that in $Q\times \cc$ 
 \begin{equation}\label{abscons}
 S_{Q,Q_1}\leq T_{Q,Q_1}\leq T,~S_{Q,Q_1}\leq S_Q\leq S\text{, and }
\ms(S_{Q,Q_1})\geq m - \dft(Q)-2\dft(Q_1) -2\e.
\end{equation}
 Let us see first why this implies the proposition: the current $S$ being uniformly laminar,
  we can use the holonomy to extend
  the current $S_{Q,Q_1}$, which is originally defined in $Q\times \cc$, to 
 $Q_1\times \cc$, and get a current we still denote by $S_{Q,Q_1}$, subordinate to both
 $S$ and $T$ in $Q_1\times \cc$ and  satisfying the last estimate in 
 (\ref{abscons}). So we get as before the following estimate in $U\times \cc$
 $$\left \langle \sum_{Q_1\in\qq_k} S_{Q,Q_1},\mathbf{1}_{U\times\cc}
 ~ \pi^* idz\wedge d\bar z
\right\rangle \geq \langle S, \mathbf{1}_{U\times\cc} ~ 
\pi^* idz\wedge d\bar z\rangle -4\e,$$ that is, $\sum_{Q_1\in\qq_k} S_{Q,Q_1}$ increases 
to $S$. On the other hand, $\sum_{Q_1\in\qq_k} S_{Q,Q_1}\leq T$ and we conclude that 
$S\leq T$ in $U\times \cc$.\\

It remains to prove our claim. The data are
 $$S_Q\leq T_Q,~\ms(S_Q)\geq m-\e\text{ and }T_{Q,Q_1}\leq T_Q,~\ms(T_{Q,Q_1})   
\geq 1-\dft(Q)-2\dft(Q_1)-\e,$$ all  these currents being uniformly laminar in $Q\times \cc$. 
Fix a global transversal $\set{c}\times \cc$, $c\in Q$, and consider the respective slices 
$m_{S_Q}$, $m_{T_Q}$, $m_{T_{Q,Q_1}}$ of  $S_Q$, $T_Q$ and $T_{Q,Q_1}$.
By the Radon-Nikodym Theorem there exists a function $0\leq f_{S_Q}\leq 1$ (resp.
$0\leq f_{T_{Q,Q_1}}\leq 1$) such that $m_{S_Q}=f_{S_Q} m_{T_Q}$ (resp. 
$m_{T_{Q,Q_1}}= f_{T_{Q,Q_1}}m_{T_Q}$). 

Let $f=\inf(f_{S_Q}, f_{T_{Q,Q_1}})$, then one has the estimate 
 $$\int f dm_{T_Q} \geq m - \dft(Q)-\dft(Q_1)-2\e.$$ Define 
 $S_{Q,Q_1}$ as the uniformly laminar current in $Q\times \cc$ subordinate to $S$, and having transverse measure  $fdm_{T_Q}$ in $\set{c}\times \cc$. 
 $S_{Q,Q_1}$ has the required properties (\ref{abscons}).\hfill $\square$\\
 
 \begin{rmk}\normalfont\label{rmk_web}
 The definition of laminar currents may be relaxed to let the disks intersect. One obtains 
 the class of {\sl web-laminar} currents, considered by Dinh \cite{dinh}, which seems to be of interest. For instance, the cluster values of a sequence of curves in $\pd$ with degree $d_n$ 
 and geometric genus $O(d_n)$ and no assumption on the singularities 
  are of this form --such a statement may easily be extracted from \cite{lamin}, and is explicit in \cite{dinh}. Moreover such currents are strongly approximable in the sense 
 that estimate (\ref{eq_goodcp}) holds, with the $T_{Q,n}$ being sums of intersecting graphs (web-uniformly laminar currents).
 
One may then define disks subordinate to a web-laminar current, and prove an analytic
 continuation theorem in the strongly approximable case,
 in the same way as above. 
 \end{rmk}
 
 \begin{rmk}\normalfont\label{rmk_relax}
 The estimate (\ref{eq_goodcp}) plays of course an important role in this section. 
 However a careful reading of proposition \ref{prop_percol} shows it can be relaxed
 by replacing $O(r^2)$ by $O(r^{1+\e})$; in particular the analytic continuation 
 statement holds in this case. 
 \end{rmk}
  
%%%%%%%%%%%%%%%%%%%%%%%%%%%%%%%%%%%%%%%%%%
%%%%%%%%%%%%%%%%%%%%%%%%%%%%%%%%%%%%%%%%%%

\section{Non self intersection}\label{sec_isect}

In this section we prove the second part of theorem \ref{thm_main}, which asserts that disks subordinate to a diffuse strongly approximable $T$ are compatible. Due to our definition \ref{def_sub} of disks subordinate to a laminar current, this is equivalent to saying that uniformly laminar currents subordinate to $T$ do not intersect non trivially. 

In contrast to the preceding section, the result here is purely local, and only uses the mass estimate of proposition \ref{prop_mass}. We first recall the statement.

\begin{thm}\label{thm_isect}
Let $T$ be a strongly approximable and diffuse laminar current in $\om\subset\cd$. Then two disks subordinate to $T$ are compatible, i.e. their intersection is either empty or open in the disk topology.
\end{thm}

A few comments are in order here. First there are simple examples of laminar currents with intersecting subordinate disks, given by sums of uniformly laminar currents with transversals of zero area (see e.g. \cite[exemple 2.2]{isect}). Thus such currents cannot be strongly approximable.

On the other hand if 
one weakens the definition of disks subordinate to $T$ to the following
 ``a disk is subordinate to $T$ if it is the union of disks appearing in the laminar representation (\ref{eq_rep}), up to a set of zero measure", then disks subordinate (in this sense)
 to a strongly approximable $T$ may intersect. For example a pencil of lines with any transverse measure satisfies (\ref{eq_mass}). We nevertheless believe this is not
 a workable definition of disks subordinate to $T$. \\ 
  
In case $T$ has continuous potential, the theorem is a consequence of 
the results of \cite{isect}. Indeed we proved that $T\wedge T=0$ in this case, so if $S_1$ and $S_2$ are uniformly laminar currents subordinate to $T$ in $\om'\subset\om$, one has $S_1\wedge S_2=0$ (currents dominated by $T$ also have continuous potential). One interesting point here is that no potential is involved; the result may thus appear as a ``geometric version" of the equation $T\wedge T=0$.\\

\begin{pf} the proof is by contradiction. So assume that
$S_1$ and $S_2$ are uniformly laminar currents in $\om'\subset\om$, with non trivial intersection, 
and such that $S_i\leq T$. 
It is no loss of generality to assume $\om'=\om$. Most intersections between the leaves of the associated laminations $\el(S_1)$ and $\el(S_2)$
are transverse by \cite[Lemma 6.4]{bls}, so 
focusing on a neighborhood of such  a transverse intersection point and reducing $\om$, 
$S_1$ and $S_2$ if necessary, we 
assume $S_1$ and $S_2$ are made up of almost parallel disks and
 that any leaf of $S_1$ is a global transversal to $\el (S_2)$.\\

Next, recall that $T$ is the increasing limit of  sums of uniformly laminar currents in 
cubes $\sum T_Q$ given by proposition \ref{prop_mass}.
The approximation  is increasing, 
so if a disk subordinate to, say, $S_1$ appears at some stage of the approximation, it will persist in all finer 
subdivisions. Moreover, the approximating currents are uniformly laminar, so in the approximation, disks do not ever intersect non trivially. 

In what follows
 the notation $\qq$ denotes  subdivisions by families of affine {\sl cubes} in $\om$, as given by  proposition \ref{prop_mass}. 
Recall also from this proposition that subdivisions may be translated since only the projections $\pi_1$ and $\pi_2$ are fixed, and estimate (\ref{eq_mass}) still holds.

There are two mutually disjoint cases. 
\begin{itemize}
\item[-] Either  at some stage of the approximation, 
one obtains a current $T_Q$, with
 $S'_1\leq T_Q\leq S_1$, such that $\supp(S'_1)\cap\supp(S_2)
\neq\emptyset$ --the case where 1 and 2 are swapped is similar. 
In this case, the disks subordinate to $S'_1$ persist in finer subdivisions, 
and the corresponding intersecting disks subordinate to $S_2$ never appear. 
\item[-] Or such a current never appears.
\end{itemize}
In both cases, some disks subordinate to $S_2$ will never appear in the approximation
process. More precisely, these correspond to the set of disks in subdivisions by cubes 
of size $r$, subordinate to $S_2$, and intersecting some fixed $S'_1\leq S_1$. We wish 
to prove that
this contradicts estimate (\ref{eq_mass}). Without loss of generality, we put 
$S'_1=S_1$, we also renormalize the transverse measures so that the measure induced 
by $S_2$ on the leaves of $S_1$ is of mass 1, and make an affine transformation
so that the projections $\pi_1$ and $\pi_2$ become orthogonal.\\

For a given subdivision by affine cubes $\qq$, and $Q\in \qq$, we denote by $\frac{Q}{2}$ 
the image of $Q$ by scaling of factor $1/2$ with respect to its center. 

\begin{lem} For every $r>0$, there exists a subdivision $\qq$ by cubes of size $r$, and $N(r)$ leaves $(L_i)_{i=1}^{N(r)}$ of $\el(S_1)$, with mutual distance $\geq 5r$, such that 
if $m_i=S_2\wedge [L_i]$ denotes the transverse measure induced by $S_2$ on $L_i$, one has 
$$\left( \sum_{i=1}^{N(r)} m_i \right) \left( \bigcup_{Q\in \qq} \frac{Q}{2}\right) \geq \frac{N(r)}{32} \underset{r\cv 0}{\longrightarrow} \infty.$$
\end{lem}

Let us see first why the lemma implies the theorem. By the reductions made so far, we know that no disk of $T_Q$ traced on a 
leaf of $\el(S_2)$ intersects the $N(r)$ leaves $L_i$. 
Moreover, for every cube $Q$ of size 
$r$, there exists a constant $c$ such that any subvariety of $Q$ intersecting $\frac{Q}{2}$ has area at least $cr^2$ (Lelong's Theorem). Thus if $L$ is any leaf of $\el(S_1)$, the 
total mass of the uniformly laminar current subordinate to $S_2$, made up of the disks 
through $\frac{Q}{2}$ is at least $(S_2\wedge[L])(\frac{Q}{2}) cr^2$. 
Since leaves at $5r$ distance cannot hit the same (affine) cube of size $r$, the preceding lemma provides us with a sum of uniformly laminar currents, subordinate to $S_2$, with 
mass greater than $\frac{c}{32}N(r)r^2$, that will never appear in the approximation process. This is a contradiction since $N(r)\cv\infty$.\\

\noindent{\bf Proof of the lemma:} first recall that
the holonomy of $\el(S_1)$ is H{\"o}lder 
continuous, so if a transverse section 
 of $\el(S_1)$ is fixed, for appropriate constants $C$ and $\tau$,  points mutually distant of $C r^\tau$ in the transversal give rise to leaves distant of
$5r$ in $\om$. Pick $N(r)$ such points in the transversal; as $r\cv 0$, $N(r)\cv\infty$ since 
$S_1$ is diffuse. 

For the associated leaves $L_i$, let $m_i$ be the measure induced by $S_2$, and 
$m=\sum_i m_i$, which is a measure of mass $N(r)$ by the normalization done before.  
It is an easy consequence of the translation invariance of Lebesgue measure and the Fubini Theorem (see \cite[lemma 4.5]{isect}) 
that there exists a translate of $\qq$ such that the mass of $m$ concentrated in $\cup \frac{Q}{2}$ is larger than 
$$\unsur{2} \frac{{\rm volume \left(\frac {Q}{2}\right)}} {{\rm volume(Q)}} N(r) = \frac{N(r)}{32},$$
which yields  the desired conclusion.
\end{pf}

%%%%%%%%%%%%%%%%%%%%%%%%%%%%%%%%%%%
%%%%%%%%%%%%%%%%%%%%%%%%%%%%%%%%%%%
\section{Measured laminations}\label{sec_su}

In this section we will reinterpret  the preceding results in a more geometric fashion, 
by constructing a {\sl  weak measured lamination} associated to a strongly approximable
current in $\pd$. This has the advantage of clarifying the question of {\sl representation} 
of laminar currents, since the measured lamination is the ``largest" 
possible representation. 
We emphasize that this construction is more
generally valid for currents satisfying the conclusions of Theorem \ref{thm_main}.
In analogy with Cantat \cite{ca}, we define those as the {\sl strongly laminar
currents}. Theorem \ref{thm_main} paraphrases then as ``strongly approximable
currents are strongly laminar".

 Next we relate closed currents subordinate to $T$ and invariant transverse
measures on its associated weak lamination (Theorem \ref{thm_sub}) --this result
really needs the mass estimate (\ref{eq_mass}).

\subsection*{Weak laminations.}
We first define a notion of weak lamination adapted to our setting. The definitions
are {\sl ad hoc} so we assume the ambient space is a 2 dimensional complex 
manifold. We fix a diffuse strongly laminar current $T$. 

\begin{defi}\label{def_flowbox}
A flow box for $T$ is the (closed) support of  a lamination $\el$ embedded in $U\simeq\dd^2$, such that in this coordinate chart $\el$ is 
biholomorphic to a lamination by graphs over the unit disk, 
and moreover satisfying  $T\rest{\el}>0$ and $\supp(T\rest{\el})=\el$.

The regular set $\mathcal{R}$ is the union of the disks subordinate to $T$, or
equivalently the union of flow boxes.
\end{defi}

The condition on the support of $T\rest{\el}$ insures that we do not consider  disks not
subordinate to $T$.  By definition, for any laminar current, the
regular set has full 
measure in $T$. 

\begin{defi}\label{def_weak}
Two flow boxes are said to be compatible if the associated disks intersect in  a compatible 
way. 

A weak lamination is a  union of a family of compatible flow boxes. We say it is 
$\sigma$-compact if there are countably many flow boxes.
\end{defi}

It turns out that this definition fits well with the theory of laminar currents, where 
no transverse topology is {\sl a priori} involved --see however the density topology below. 
After this paper was written, we realized that a similar definition already appears
in Zimmer \cite{z}.
Given a weak lamination, one easily define leaves, 
 as in the usual case; if the weak lamination is $\sigma$-compact then the leaves are 
 $\sigma$-compact (see e.g. \cite{cc} or \cite{these}).
 
 We say a closed set $\tau\ni x$ is a {\sl local transversal} to the weak lamination 
 at $x$ if it is 
 a local transverse section  {\sl in a  flow box}. Due to the compatibility condition, 
 this is independent of the choice of the flow box containing $x$. One then defines 
 holonomy maps between transversals as in the usual case, and one may speak 
 of holonomy invariant transverse measures, that is, a collection of measures on all
 transversals, invariant by holonomy (see Sullivan \cite{su}, Ghys \cite{ghys}). 

The following proposition is a reformulation of Theorem \ref{thm_main}.

\begin{prop}\label{prop_flowbox}
Let $T$ be a laminar current satisfying the conclusions of Theorem \ref{thm_main}
-- a strongly laminar current. 
Then the regular set $\mathcal{R}$ has the structure of a  weak lamination, and
$T$ induces a holonomy  invariant transverse measure on $\mathcal{R}$.
\end{prop}

\begin{pf} the only non trivial statement is the existence of the invariant transverse 
measure. Note first that $T$ induces an invariant transverse measure on each flow box
$\el$, since $T\rest{\el}$ is uniformly laminar. Now if two flow boxes $\el_1$ and $\el_2$
 have non trivial intersection --compatible due to the 
 non intersection of disks subordinate to $T$-- , the transverse 
measures coincide on common transversals: just construct a flow 
box $\el$ from this transversal, subordinate to both $\el_1$ and $\el_2$, and apply the 
analytic continuation theorem again.
\end{pf}

Notice that this result gives a {\sl natural} representation of $T$ as an integral 
over families of disks, since the definition of $\mathcal{R}$ does not involve
any choice: we take all disks subordinate to $T$. 
 This means the class of 
strongly laminar currents should be a reasonable intermediate class between general 
and uniformly laminar currents.\\

The following intuitive  proposition asserts the transverse mass of a flow box is 
computed using wedge products. We know that if $\el$ is a flow box,
$T\rest{\el}$ is uniformly laminar, so if $\tau$ is any  global
holomorphic transversal, 
the transverse mass of $\el$ is given by 
$m=\m(T\rest{\el}\wedge [\tau])$, which is easily proved to be a well defined wedge 
product. The expected thing is that $m=(T\wedge [\tau])(\el\cap\tau)$ provided the
wedge product is well defined, which is almost true: this is the
content of the next proposition. Note that using the techniques of
section \ref{sec_pluripotential} one may replace the smooth uniformly
laminar currents in the proof  by
uniformly laminar currents not charging  pluripolar sets.   

\begin{prop} \label{prop_wedge}
Let $(\tau_\lambda)_{\lambda\in\dd}$ be a smooth family of disjoint global 
holomorphic transversals
to $\el$. Then for almost every $\lambda\in\dd$, $T\wedge[\tau_\lambda]$ is 
a well defined positive measure
and the transverse mass of $\el$ is $(T\wedge[\tau_\lambda])(\el\cap\tau_\lambda)$.
\end{prop}

\begin{pf} let $\psi$ be any positive test function in $\dd$ and let $S$ be the smooth uniformly laminar current $S=\int_\dd [\tau_\lambda]\psi(\lambda)d\lambda.$ Since $S$
is smooth, the wedge product $S\wedge T$ is well defined and described by the
geometric intersection of disks constituting $S$ and $T$ (for more details on this topic
see \cite{isect}), i.e. there is a laminar representation of $T$,
$T= \int_{\mathcal{A}}[D_a]d\mu(a)$ such that 
$$S\wedge T = \int_{\mathcal{A}\times\dd} [\tau_\lambda\cap D_a]\psi(\lambda)
d\lambda d\mu(a).$$
Indeed since $S$ is 
smooth one has 
$$T\wedge S= \int_\mathcal{A} ( [D_a]\wedge S) d\mu(a) = \int_\mathcal{A} 
\left(\int_{D_a} S\right) d\mu(a),$$
and $S$ being uniformly laminar, $[D_a]\wedge S$ is a geometric intersection. 

It is a  classical fact
that for a.e. $\lambda$, the wedge product $T\wedge [\tau_\lambda]$ is well defined,
and by the preceding argument it is geometric. We conclude by using the fact that disks
subordinate to $T$ are compatible, so through every point in $\el\cap\tau_\lambda$, 
the only disk subordinate to $T$ is the corresponding leaf of $\el$.
\end{pf}

\subsection*{Su's construction.}
In the specific case of strongly approximable currents in $\pd$ one has a little bit more 
information, since the slices by generic lines give probability measures, yielding the
notion of defect. These ``reference" measures  allow one,
following Meiyu Su \cite{su}, to produce a topology --the 
{\sl density topology}-- in which $\mathcal{R}$ becomes a genuine lamination. 
This actually does not give more structure on the weak lamination, since the topology
is canonically associated to the measurable structure. 
We do not give full details, the reader is referred to
\cite{su} and \cite{these} instead.\\
 
 We fix a diffuse 
strongly approximable $T$ in $\pd$, and a linear projection $\pi$ such that
the condition on projections described in 
 definition \ref{def_sa} holds  for $\pi$ and
 the set of vertical disks has zero measure. So we can define the defect function as before, and the slice $m^z$ of $T$ by the fiber $\pi^{-1}(z)$ for every  fiber of 
 zero defect (by an increasing limit process; note that in general slicing is only
 defined for fibers  outside a polar set). By definition of the defect, the regular set 
 $\mathcal{R}$ has full transverse measure in regular fibers. 
 
 Now pick a regular fiber $\pi^{-1}(z)$, and $A$ a measurable subset of $\pi^{-1}(z)$. 
 Recall $w$ is a density point of $A$, relative to $m^z$, iff
 $$\frac{m^z (A\cap B(w,r))}{m^z(B(w,r))}\underset{r\cv 1}\longrightarrow   1; $$
Lebesgue's Theorem asserts that  if $m^z(A)>0$, almost every point in $A$ is a density 
point; it holds in the case of Radon measures in Euclidean space, see \cite{mat}.

 We can now define the {\sl density topology} by specifying its open sets.

\begin{defi}
A subset $A \in \pi^{-1}(z)$ is d-open if it is empty  or
if $A$ is measurable, $m^z(A)>0$, and every $w\in A$ is a density point.
\end{defi}

One easily checks this defines  a topology in $\pi^{-1}(z)$. Given a  flow box made 
of graphs for the projection $\pi$, we define the density topology on the flow box
as the product of the usual topology along the leaves and the (restriction of the)
density topology on a given vertical transversal. Since  holonomy   maps are 
continuous and preserve the measures $m^z$ restricted to the flow box
(thus preserve density points), this is 
independent of the transversal chosen. Moreover the density topologies on 
intersecting flow boxes coincide by compatibility and invariance  of the transverse measure. \\

We collect the following simple facts pertaining to  the d-topology:
\begin{itemize}
\item[-] sets of measure zero are d-closed sets of empty interior. 
In particular removing (possibly countably) many sets of zero measure does not affect
 the d-open property.
 \item[-] open sets of positive measure are d-open, so that the d-topology refines 
 the ambient topology in flow boxes.
 \item[-] the d-topologies induced by distinct generic linear projections $\pi$
 coincide.
 \end{itemize}
 
 We can now formulate Su's Theorem in our setting.
 
 \begin{thm}\label{thm_su}
There exists a lamination $\el$ and an injection $i:\el\hookrightarrow \supp T$ 
continuous along the leaves and respecting the laminar structure, with image of full
measure, and full transverse measure on each transversal. 

The lamination $\el$ has an invariant transverse measure, and if $T(\el)$ denotes 
the associated foliated cycle, one has $i_*(T(\el))=T$.
\end{thm}  
 
\begin{pf} we construct a d-open subset $\mathcal{R'}\subset\mathcal{R}$ of full 
transverse measure, and saturated with respect to  the weak lamination on $\mathcal{R}$.
Let $(B_m)_{m\geq 0}$ be a covering of $\mathcal{R}$ by flow boxes of
 positive measure. The d-interior of $B_m$ is denoted by $\dint{B_m}$. For each box 
 $B_k$, one removes from $\dint{B_k}$ all the plaques corresponding to leaves containing plaques 
 of $\cup_{m\geq 0} B_m\backslash \dint{B_m}$ ($k$ itself is included in the union 
 since a leaf may intersect $B_k$ several times). The set of removed  plaques has zero 
 transverse measure so it remains a d-open subset $B_k'\subset
 B_k$. Let $\mathcal{R}'=\cup B'_k$, which is d-open and saturated by construction. \\
 
 Now for every $x\in\mathcal{R}'$, $x\in B_k'$ for some $k$, and $B_k'$ is a foliated
 d-neighborhood of $x$, so $\mathcal{R}'$ supports a natural lamination $\el$, which has the desired properties.
 \end{pf}
 
 \subsection*{Subordinate transverse measures} In this paragraph we relate 
subordination at the levels of transverse measures and closed laminar currents. 

\begin{thm}\label{thm_sub}
Let $T$ be a strongly approximable current on $\pd$, and $\mu$ be the 
induced transverse measure on its associated weak lamination $\mathcal{R}$. 
Then to every invariant transverse measure $\mu'\leq \mu$ on $\mathcal{R}$, 
the foliation cycle induced by $\mu'$ on $\mathcal{R}$
 corresponds to a closed strongly laminar current $T'\leq T$ in $\pd$.
\end{thm}

\begin{cor}
If $T$ is an extremal current, the transverse measure $\mu$ is ergodic.
\end{cor}

Recall  that ``ergodic" means every 
saturated set $\mathcal{R}'\subset\mathcal{R}$ has either zero or full measure. 
This corollary has non trivial dynamical consequences that 
will be developed in further work.
The theorem means that there is a good (one way)
correspondence between closed positive 
currents on the weak lamination $\mathcal{R}$ and closed positive currents on the ambient space. We do conjecture the converse also holds, that is, {\it every closed 
positive current $T'\leq T$ is the foliation cycle of some invariant transverse measure 
$\mu'\leq \mu$.} This would imply for instance that $T$ is extremal if and only if
its transverse measure is ergodic. This conjecture seems to raise rather delicate 
problems of analytic type. \\

\begin{pf} the result is local in some open $\om$.
Pick two linear projections $\pi_1$ and $\pi_2$ satisfying definition
\ref{def_sa}, and for each basis of projection, consider 3 overlapping subdivisions 
by squares of size $r$, as in section \ref{sec_cont}. We then form 9 overlapping
subdivisions $\qq_1,\ldots,\qq_9$ by affine cubes 
of  size $r$ from the projections and the squares, as in 
proposition \ref{prop_mass}. For each $\qq_i$, one gets a current $T_{\qq_i}$, 
uniformly laminar in each cube, and such that $\m (T-T_{\qq_i})=O(r^2)$.\\

We are given a measure $\mu'$ on the weak lamination $\mathcal{R}$, with 
$\mu'\leq \mu$. This  means that for each local transversal $\tau$ (in a
flow box by definition of the transversals), there exists a function $f_\tau$, 
$0\leq f_\tau\leq 1$, with $\mu'=f_\tau\mu$. Since $\mu'$ is invariant by holonomy,
the functions $f_\tau$ patch together 
to a global $0\leq f\leq 1$ on $\mathcal{R}$, 
constant along the leaves. The associated foliated cycle is the current $fT$. We have to
 prove it is closed.\\
 
So let $\phi$ be a test 1-form in $\om$. We consider a partition of unity 
$(\theta_Q)_{Q\in\qq}$ subordinate 
to the covering $\qq=
\qq_1\cup \cdots \cup \qq_9$ of $\om$. It is easily seen that since the 
cubes have size $r$, the functions $\theta_Q$ may be chosen to have derivatives 
uniformly bounded by $C/r$, with $C$ independent of $r$. Now,
$$\bigl\langle fT, d\phi\bigl\rangle 
= \bigl\langle fT , d\bigl( \sum_{Q\in\qq} \theta_Q \phi \bigl) 
\bigl\rangle= \sum_{i=1}^9 \bigl\langle fT , d\bigl( \sum_{Q\in\qq_i} \theta_Q \phi \bigl) 
\bigl\rangle, $$  and for each of the nine terms of the sum, we replace $fT$ by
$fT_{\qq_i} + f(T-T_{\qq_i})$. The important fact is that since $f$ is constant along the 
leaves, for each $Q\in\qq_i$ the current $fT_{\qq_i}\rest{Q}$ is closed in $Q$, so we get
$$\bigl\langle fT_{\qq_i} , d\bigl( \sum_{Q\in\qq_i} \theta_Q \phi \bigl) 
\bigl\rangle = \sum_{Q\in\qq_i} \bigl\langle fT_{\qq_i}\rest{Q}
, d(\theta_Q \phi) \bigl\rangle = 0,$$
since  $\theta_Q$ has compact support in $Q$. On the other hand 
$$\abs{\bigl\langle f(T-T_{\qq_i}) , d\bigl( \sum_{Q\in\qq_i} \theta_Q \phi \bigl)}\leq 
\m\bigl(f (T-T_{\qq_i})\bigl) \sup_{Q\in\qq_i} \norm{ d(\theta_Q \phi)}\leq 
\m\bigl(T-T_{\qq_i}\bigl) O(\unsur{r}) = O(r).$$
This implies $\langle fT, d\phi\rangle=0$ and the theorem follows. 
\end{pf}
 
%%%%%%%%%%%%%%%%%%%%%%%%%%%%%%%%%%%%
%%%%%%%%%%%%%%%%%%%%%%%%%%%%%%%%%%%%
 
\section{Pluripotential theory and laminar currents}\label{sec_pluripotential}

We give in this paragraph a few applications of the foregoing study. More precisely we first prove that the potential of a strongly approximable 
laminar current $T$ is either harmonic or identically $-\infty$ on
almost all disks subordinate 
 to $T$ (leaves of the induced measured lamination). We also exhibit a decomposition
 of a strongly approximable laminar current into a sum of two closed laminar currents, 
 one essentially supported on a pluripolar set and the other not charging pluripolar sets.
 
\subsection*{Some results on uniformly laminar currents}
We first collect some useful results on uniformly laminar currents. 
The proofs only use a few simple ideas from 
1-variable classical potential theory. Our first goal
is the following proposition, although the intermediate lemmas may be of
independent interest.

\begin{prop}\label{prop_ul}
Let $T$ be a uniformly laminar current, given as the 
integral of holomorphic graphs
in the bidisk, $T = \int [\Gamma_\alpha] d\mu(\alpha)$. Assume $T$
does not give mass to pluripolar sets. Then  $T$ can be written as a countable
sum $T=\sum
T_j$, where the $T_j =  \int [\Gamma_\alpha] d\mu_j(\alpha)$ have
continuous potential and disjoint support. 
\end{prop}

The main ingredient of the proof is the following 1-variable result,
which may be found for instance in H{\"o}rmander's book 
\cite[Theorem 3.4.7]{h}.

\begin{prop}\label{prop_1V}
Let $\mu$ be a positive measure with compact support in $\cc$. Assume
the logarithmic potential $G_\mu(z)=\int\log\abs{z-\zeta}d\mu(\zeta)$ 
satisfies $G_\mu>-\infty$
$\mu$-a.e.--this is true in particular if $\mu$ does not charge polar
sets. Then there exists a sequence of disjoint compact subsets $(K_j)$
such that  $\mu=\sum \mu\rest{K_j}$ and for each $j$, $\mu\rest{K_j}$
has continuous potential.
\end{prop} 

The proposition is a consequence of Lusin's theorem, together with the
 ``continuity principle'' for logarithmic potential.

We proceed, in several steps, to
the proof of proposition \ref{prop_ul}. We denote by $\el$ the
 lamination by
horizontal graphs in the bidisk, associated to $T$. 
The family of (vertical) disks
$\set{z}\times \dd$ is a family of global transversals to the
lamination. Let $h^z=h^{0,z}$  be the holonomy map from $\el\cap
(\set{0}\times \dd)$ to $\el\cap
(\set{z}\times \dd)$, and similarly $h^{z,z'}$. 
We identify an abstract transverse measure $\mu$ on
$\el$ with its image in $\set{0}\times \dd$, so that the parameter
$\alpha$ is identified with the point $(0,\alpha)$, and let $\mu^z=
(h^z)_*\mu$ be the push forward of $\mu$ in $\set{z}\times \dd$. 

\begin{lem}\label{lem_potential}
Let $T = \int [\Gamma_\alpha] d\mu(\alpha)$ as above. Then the
function $$u_T: (z,w)\longmapsto \int_{\set{z}\times\dd} \log \abs{w-\zeta}
d\mu^z(\zeta)$$ is a plurisubharmonic potential for $T$.
\end{lem}
\begin{pf} classical, we include it for completeness. Let
$w=\varphi_\alpha(z)$ be the equation of the graph
$\Gamma_\alpha$. Then 
\begin{equation}\label{eq_pot}
u(z,w)= \int \log\abs{w-\varphi_\alpha(z)} d\mu(\alpha)
\end{equation} is a
potential for $T$. Now the holonomy map $h^z$ maps $(0,\alpha)$ to
$(z,\varphi_\alpha(z))$, so for any continuous function $F$
on $\set{z}\times\dd$
$$ \int_{\set{z}\times\dd} F(\alpha) d\mu^z(\alpha)= 
 \int_{\set{0}\times\dd} F(\varphi_\alpha(z)) d\mu(\alpha)$$ and 
writing $\log\abs{w-\cdot}$ as a decreasing sequence of continuous
functions, we get $u=u_T$.
\end{pf}

There is a good correspondence between continuity properties 
of the potentials of the transversal measures and the current itself:
this is the content of the next lemma.

\begin{lem} \label{lem_cont_holonomy}
Assume $\mu$ has continuous potential as a measure on
$\set{0}\times\dd$. Then for every $z$, $\mu^z$ has continuous potential.
Moreover the above defined potential $u_T$ is continuous.
\end{lem} 

\begin{pf} the first assertion is a consequence of the following: 
the class of plane measures with 
continuous potentials is preserved by bi-H{\"o}lder continuous
homeomorphisms. Indeed let $\mu$ be a plane positive measure with compact
support, and 
$$k_\mu(z,r) = \abs {\int_{B(z,r)} \log\abs{z-\zeta}d\mu(\zeta)}.$$
Let $c\geq 2$; using the fact that $B(w,c\abs{z-w})\subset
B(w,(c+1)\abs{z-w})$ and the mean value inequality for the logarithm 
one easily gets 
$$\abs{G_\mu(z)-G_\mu(w)} \leq \unsur{c} + k(z,c\abs{z-w}) + 
k(w,(c+1)\abs{z-w}).$$ Taking for example $c=\abs{z-w}^{-1/2}$,
one deduces the following result: if $k(z,r)\cv 0$ locally uniformly 
in $z$   as $r\cv 0$, then $G_\mu$ is continuous. Using similar
estimates it is proven by
Shvedov \cite{shvedov} that the converse is also true. \\

Assume now $\mu$ has  continuous potential in $\set{0}\times
\dd$. Then $k_\mu(x,r)\cv 0$ uniformly by the Shvedov result; moreover 
the holonomy map $h^z$ associated to $\el$ is 
H{\"o}lder continuous, as well as its inverse, 
say of exponent $\alpha$, and we get 
$k_{\mu^z}(w,r)\leq Ck_\mu((h^z)^{-1}(w),Cr^\alpha)$. This implies
$\mu^z$ has continuous potential also --note that the modulus of
continuity is uniform.\\

It remains to prove continuity of $u_T$ as a function of $(z,w)$. First, we
extend the lamination $\el$ to a neighborhood of $\supp(T)$ using
Slodkowski's theorem. Now it follows from formula (\ref{eq_pot})  that the
potential $u_T$ is harmonic or identically  $-\infty$ along the leaves.
Under the hypothesis of the theorem we know that  
the restrictions to the  slices are continuous and using the H{\"o}lder
property of holonomy again, one easily gets that
$u_T$ is bounded.  
Then we split
$$u_T(z,w)  - u_T(z',w')= \left(u_T(z,w)- u_T(z, h^{z',z}(w')\right) + 
\left(  u_T(z, h^{z',z}(w') -  u_T(z',w')\right),$$
where the first term on the right hand side is small because of
the continuity  of $\zeta \mapsto u_T(z,\zeta)$, and the second
because $z\mapsto  u_T(z, h^{z',z}(w'))$ is a uniformly bounded 
harmonic function, hence uniformly Lipschitz.
\end{pf} 

Recall that $X\subset\cc$ (resp. $X\subset\cd$) 
is {\sl polar} (resp. {\sl pluripolar})
if $X\subset \set{u=-\infty}$ where $u$ is a subharmonic
(resp. plurisubharmonic) function, not identically equal to $-\infty$. 

\begin{lem}\label{lem_polar}
Let $X$ be a subset of $\set{0}\times \dd$, and $\widehat{X}$
the set 
saturated from $X$ by the lamination $\el$ (i.e. the set of leaves through
$X$). Then $X$ is polar iff $\widehat{X}$ is pluripolar.
\end{lem}
 
\begin{pf} note first that the holonomy map preserves the class of closed
polar subsets of the fibers $\set{z}\times \dd$. A way to prove this
is to use the following characterization of polar sets  (transfinite
diameter zero, see Tsuji \cite{t}) : $X$ is polar iff 
$$\lim_{n\cv\infty} ~\sup\set{ \prod
\abs{x_i-x_j}^{2/n(n-1)},~x_1,\ldots,x_n \in X }=0 $$ 
This condition is stable under bi-H{\"o}lder homeomorphisms. Another method is to
use lemma \ref{lem_cont_holonomy} and the fact that a non polar
compact set
carries a measure with continuous potential \cite[Theorem
3.4.5]{h}.

If $X$ is not closed (polar sets are $G_\delta$ sets in general), use the fact that 
$X$ is polar iff for every compact $K\subset X$, $K$ is polar, and rather
transport the compact subsets. \\

Now assume $\widehat{X}$ is pluripolar. Then
$\widehat{X}\subset\set{u=-\infty}$ for some non degenerate
p.s.h. function in $\dd^2$. Hence for almost every slice
$\set{z}\times\dd$, $u\rest{\set{z}\times\dd}\not\equiv -\infty$ and
$\big(u\rest{\set{z}\times \dd}\big)\big( \widehat{X}\cap (\set{z}\times\dd)\big)
\equiv -\infty$, so
$\widehat{X}\cap (\set{z}\times\dd)$ is polar. The preceding
observation implies $X=\widehat{X}\cap (\set{0}\times\dd)$ is polar.\\

Conversely, assume $X$ is polar. Then by \cite[theorem 3.4.2]{h}
there exists a positive measure $\mu$ supported on $X$ such that $X\subset\set
{G_\mu=-\infty}$. Consider the following plurisubharmonic function
$$u(z,w)=\int_{\set{z}\times\dd} \log \abs{w-\zeta}
d\mu^z(\zeta).$$ On each leaf, $u$ is harmonic or identically
$-\infty$. We thus get  $\widehat{X}\subset \set{u=-\infty}$.
\end{pf} 

From these lemmas one easily deduces the proof of proposition
\ref{prop_ul}. Assume $T$ does not charge pluripolar sets. Then the
transverse measure does not charge polar subsets of 
$\set{0}\times\dd$ by lemma \ref{lem_polar}. 
Write $\mu=\sum \mu_j$ as given by proposition
\ref{prop_1V}, and $T=\sum T_j$ according to this decomposition. By lemma
\ref{lem_cont_holonomy}, $T_j$ has continuous potential. \hfill $\square$\\    
 
 The next proposition gives a decomposition of a uniformly laminar current in the bidisk
  as a sum of two parts, one giving mass to pluripolar sets, the other not. It will be used in Theorem \ref{thm_decomp}.
 
 \begin{prop}\label{prop_decomp_ul}
 Let $T$ be a uniformly laminar current, integral of holomorphic graphs
in the bidisk, $T = \int [\Gamma_\alpha] d\mu(\alpha)$. Then $T$ admits 
a unique decomposition  as a sum 
$T=T'+T''$ of uniformly laminar currents, with $T'$ not charging pluripolar sets, and $T''$ 
giving full mass to a pluripolar set.
\end{prop}
 
 \begin{pf} uniqueness is obvious: if $T=T_1'+T_1''= T_2'+T_2''$, just write 
 $T_1'-T_2'=T''_2-T''_1$. 
  We first decompose the transverse measure, and then apply the preceding
 lemmas. Let $\mu^0$ be the slice of $T$ by $\set{0}\times \dd$, as before, and $u_0
 =u_T(0,\cdot)$ be the logarithmic potential of $\mu^0$. Then 
 $$\mu^0=\mu'+\mu''= \mu^0\rest{\set{u_0>-\infty}}+\mu^0\rest{\set{u_0=-\infty}}.$$  
 Let $v$ be the logarithmic potential of $\mu'$. Since $\mu'\leq\mu^0$, $u_0-v$ is 
subharmonic, so $v\geq u_0 +O(1)$, and $v$ is finite $\mu'$-a.e. By proposition \ref{prop_1V} above, $\mu'$ does not charge polar sets; moreover $\mu'$ is  a sum of measures with continuous potential. On the other hand $\mu''$ gives full mass to the 
polar set $\set{u_0=-\infty}$.

Now  decompose $T=T'+T''$ according to this decomposition of $\mu^0$. By lemma \ref{lem_polar} above, $T''$ has full measure on a pluripolar set.
 Moreover, since $\mu'$ is a 
sum of measures with continuous potential, we get an analogous decomposition for $T'$
by lemma \ref{lem_cont_holonomy}, and $T'$ does not charge pluripolar sets.
\end{pf}

 \subsection*{The potential along the leaves} Recall that a disk $\Delta$ is subordinate to 
 $T$ if it is subordinate to a uniformly laminar $S\leq T$  in $\om'\subset\om$. Notice 
 in the following theorem that the condition of being harmonic or $-\infty$ on $\Delta$ is clearly independent of  potential chosen for $T$.
 
 \begin{thm}\label{thm_harmo} 
 Let $T=dd^cu$ be a  diffuse strongly approximable laminar current in 
 $\om$. Then for almost every disk $\Delta$  subordinate to $T$, with respect 
  to the transverse measure, 
 either $u\rest{\Delta}$ is harmonic, or $u\rest{\Delta}\equiv-\infty$.
 \end{thm}
 
 We remark that there are disks on which $u$ is harmonic, without being subordinate to the current in our sense. For example if $T$ is a current made up of 
 a measured family of disjoint branched coverings of degree 2 over the unit disk, say branched over $0$, and accumulating on the horizontal line (this is called a {\sl folded uniformly laminar current} 
 in \cite{these}), then the estimate (\ref{eq_mass}) is satisfied, 
 and the potential of $T$ is
  harmonic on the horizontal line, even if there is no laminated set of positive measure  
 containing it. 
 
 It would be interesting to understand more about the disks in $\supp(T)$ such that 
 $u\rest{\Delta}$ is harmonic. It seems that such disks should be
 ``tangent" to $T$ in some sense. 
 
 Also, we believe the result should be true  for {\sl every} disk subordinate to $T$. 
 Of course this is true if $u$ is continuous.\\
 
 \begin{pf} we have to prove that if $S\leq T$ is a uniformly laminar current, 
 for a.e. disk $\Delta$ subordinate to $S$, $u\rest{\Delta}$ is harmonic or $-\infty$. Reducing $\om'$ if necessary we may 
 assume $S$ is made up of graphs over some disk. We apply propositions 
 \ref{prop_decomp_ul} and \ref{prop_ul} to $S$ and get $S=S'+S''$; moreover 
 $S'=\sum S_j$, with the uniformly laminar currents $S_j$ of disjoint 
support and continuous potential.  \\

We proved in \cite[Remark 4.6]{isect} that if $S_j$ is uniformly laminar with continuous potential, the wedge product $S_j\wedge T$ is geometric, i.e. described by the geometric intersection of the disks constituting the current; moreover by Theorem \ref{thm_main}, disks subordinate to $T$ do not intersect. So $S_j\wedge T=0$. This means exactly that 
$u$ is harmonic on a.e. disk of $S_j$.  
 
On the other hand, we claim that $u\equiv-\infty$ on $S''$-a.e. leaf. Indeed
$S''$ gives full mass to a pluripolar set, so if $u_{S''}$ denotes
 the logarithmic potential of $S''$ as in lemma \ref{lem_potential},
 one has $u_{S''}\equiv-\infty$ on 
a.e. leaf of $S''$, for if $u_{S''}$ was finite on a set of positive measure on some 
transverse section, say $\set{z}\times \dd$, we could construct  a measure with continuous potential subordinate to $(\mu^z)''=S''\wedge[\set{z}\times\dd]$, which is 
impossible. We conclude that $u\equiv-\infty$ on almost every leaf of $S''$
because $T\geq S''$ implies $u\leq u_{S''}+ O(1)$.
 \end{pf}
 
 \subsection*{A canonical decomposition}
 The following result is reminiscent of both  the Skoda-El Mir extension Theorem 
 and Siu's decomposition Theorem for positive closed currents. It takes in our case a
 particularly complete form.
  
 \begin{thm}\label{thm_decomp}
 Let $T$ be a strongly approximable laminar current in $\pd$. Then
 there exists a  unique decomposition of $T$ as a sum of positive
 closed laminar currents $T=T'+T''$, 
 where $T'$ does not charge pluripolar sets, and
 $T''$ gives full measure to a pluripolar set. 
 Moreover $T'$ and $T''$ correspond to foliation cycles on the weak lamination 
 induced by $T$.
  \end{thm}
 
 In particular if the current $T$ is extremal, only one of $T'$ and $T''$ can appear.\\
 
 \begin{pf} note first that uniqueness is obvious. The proof actually implies $T''$ gives
 full mass to  a countable union of locally pluripolar sets, which is globally pluripolar in
 the special case of $\pd$, due to a theorem of Alexander. 
 
 The result is an easy consequence of theorem \ref{thm_sub}. Indeed we
 saw in proposition \ref{prop_decomp_ul} that a uniformly laminar current $S$ in a flow box
 admits a canonical decomposition $S=S'+S''$; this decomposition corresponds to 
 a decomposition of the transverse measure. So for each transversal $\tau$
 (as defined in
 section \ref{sec_su}), the measure  $\mu_\tau$ induced by $T$ has a decomposition
 $\mu_\tau' +\mu_\tau''$, which is  holonomy equivariant. Thus the transverse
 measure writes as $\mu=\mu'+\mu''$, and applying theorem \ref{thm_sub} gives the result.
 \end{pf}

\bigskip

\noindent{\sc \small Institut de Math{\'e}matiques de Jussieu,
Universit{\'e} Denis Diderot,
Case 7012,
2 place Jussieu,
75251 Paris cedex 05, France.}\\
{\tt \small dujardin@math.jussieu.fr}
  
\end{document}